# Consistency Spaces

Kerry M. Soileau

February 17, 2012

### ABSTRACT

We introduce the concept of a consistency space. The idea of the consistency space is motivated by the question, Given only the collection of sets of sentences which are logically consistent, is it possible to reconstruct their lattice structure?

## 1. INTRODUCTION

Let $X$ be a nonempty set. We define a <u>consistency space</u> $(X, \wp)$. $\wp$ is a nonempty collection of subsets of $X$ which satisfy the following conditions:

1. $X \notin \wp$

2. If $A \in \wp$ and $B \subseteq A$, then $B \in \wp$

Since $\wp$ is nonempty, condition 2 implies $\varnothing \in \wp$.

Each member of $\wp$ is called a <u>consistent set</u>, the remaining subsets of $X$ not in $\wp$ are called <u>inconsistent sets</u>. Notice that any superset of an inconsistent set must be inconsistent.

Next we define an equivalence relation on subsets of $X$.

We say that $A \sim B$ if for every $C \subseteq X$, $A \cup C \in \wp$ if and only if $B \cup C \in \wp$. In particular, if $A \sim B$, then $A \in \wp$ if and only if $B \in \wp$.

<u>Definition 1</u>: We say that a consistency space $(X, \wp)$ is <u>complete</u> if for every nonempty subset $A \subseteq X$, there exists an element $x_A \in X$ such that $A \sim \{x_A\}$.

<u>Proposition 1</u>: For any $C \subseteq X$, if $A \sim B$ then $(A \cup C) \sim (B \cup C)$.

Proof: Suppose $A \sim B$ and fix $C \subseteq X$. Suppose $(A \cup C) \cup D \in \wp$. Then since $A \sim B$ and $A \cup (C \cup D) \in \wp$, it follows that $B \cup (C \cup D) \in \wp$. But this just means $(B \cup C) \cup D \in \wp$, so $(A \cup C) \cup D \in \wp$ implies $(B \cup C) \cup D \in \wp$. Similarly $(B \cup C) \cup D \in \wp$ implies $(A \cup C) \cup D \in \wp$, hence $(A \cup C) \sim (B \cup C)$. ∎

Definition 2: If it exists, the underline{negation} of $x \in X$, denoted $\bar{x}$, is any element of $X$ satisfying the following:

1. $\{x, \bar{x}\} \notin \wp$

2. If $\{x, z\} \notin \wp$ then $\{\bar{x}, z\} \sim \{z\}$

3. If $\{\bar{x}, z\} \notin \wp$ then $\{x, z\} \sim \{z\}$

Proposition 2: If it exists, the negation of $x$ is unique in the sense that for any $y_1$ and $y_2$ both satisfying the criteria for $\bar{x}$, it must be true that $\{y_1\} \sim \{y_2\}$.

Proof: Fix $x$. Suppose $y_1$ and $y_2$ both satisfy the criteria for $\bar{x}$. Then by (1), $\{x, y_2\} \notin \wp$, hence by (2) we have $\{y_1, y_2\} \sim \{y_2\}$. Similarly, $\{x, y_1\} \notin \wp$, hence $\{y_2, y_1\} \sim \{y_1\}$. Finally, by transitivity of $\sim$ we conclude $\{y_1\} \sim \{y_2\}$. ∎

In the following we restrict our attention to consistency spaces $(X, \wp)$ for which $\bar{x}$ exists for every $x \in X$.

Proposition 3: $\{\overline{(\bar{x})}\} \sim \{x\}$.

Proof: First, by (1) we have $\{x, \bar{x}\} \notin \wp$ hence $\{\bar{x}, x\} \notin \wp$ so by (2) we have $\{\overline{(\bar{x})}, x\} \sim \{x\}$ and hence $\{x, \overline{(\bar{x})}\} \sim \{x\}$. Next, by (1) we have $\{\bar{x}, \overline{(\bar{x})}\} \notin \wp$ so by (3) we have $\{x, \overline{(\bar{x})}\} \sim \{\overline{(\bar{x})}\}$ and hence $\{\overline{(\bar{x})}\} \sim \{x, \overline{(\bar{x})}\}$. Finally, the transitivity of $\sim$ implies $\{\overline{(\bar{x})}\} \sim \{x\}$, as desired. ∎

Proposition 4: For any $x$ and $y$, $\{x, y, \bar{y}\} \notin \wp$.

Proof: For a contradiction, suppose $\{x, y, \bar{y}\} \in \wp$. Since $\{y, \bar{y}\} \subseteq \{x, y, \bar{y}\}$, we must have $\{y, \bar{y}\} \in \wp$. But this is absurd, hence $\{x, y, \bar{y}\} \notin \wp$. ∎

Proposition 5: For any $x$ and $y$, if $\{y, \bar{y}\} \sim \{z\}$, then $\{x, \bar{z}\} \sim \{x\}$.

Proof: We have immediately that $\{y, \bar{y}\} \notin \wp$. Since $\{y, \bar{y}\} \subseteq \{x, y, \bar{y}\}$, we must have $\{x, y, \bar{y}\} \notin \wp$. Since $\{y, \bar{y}\} \sim \{z\}$ and $\{y, \bar{y}\} \cup \{x\} \notin \wp$, we must have $\{z\} \cup \{x\} \notin \wp$, i.e. $\{z, x\} \notin \wp$. By Definition 2 we have $\{\bar{z}, x\} \sim x$, as desired. ∎

Proposition 6: $\{x\} \sim \{y\}$ if and only if $\{x, \bar{y}\} \notin \wp$ and $\{y, \bar{x}\} \notin \wp$.

Proof: ($\Rightarrow$) Suppose $\{x\} \sim \{y\}$. Then $\{x\} \cup \{\bar{y}\} \sim \{y\} \cup \{\bar{y}\}$, i.e. $\{x, \bar{y}\} \sim \{y, \bar{y}\}$. Since $\{y, \bar{y}\} \notin \wp$, it follows that $\{x, \bar{y}\} \notin \wp$. The claim that $\{y, \bar{x}\} \notin \wp$ is similarly demonstrated.

($\Leftarrow$) Suppose $\{x, \bar{y}\} \notin \wp$ and $\{y, \bar{x}\} \notin \wp$. Then $\{\bar{y}, x\} \notin \wp$, hence $\{y, x\} \sim \{x\}$, i.e. $\{x, y\} \sim \{x\}$. Also, $\{\bar{x}, y\} \notin \wp$, hence $\{x, y\} \sim \{y\}$. By transitivity of $\sim$, we have $\{x\} \sim \{y\}$, as desired. ∎

Proposition 7: $\{x\} \sim \{y\}$ if and only if $\{\bar{x}\} \sim \{\bar{y}\}$.

Proof: $(\Rightarrow)$ Suppose $\{x\} \sim \{y\}$. Then by Proposition 6, we have $\{x, \bar{y}\} \notin \wp$ and $\{y, \bar{x}\} \notin \wp$.

By Definition 2, $\{x, \bar{y}\} \notin \wp$ implies $\{\bar{x}, \bar{y}\} \sim \{\bar{y}\}$. By the same Definition, $\{y, \bar{x}\} \notin \wp$ implies $\{\bar{y}, \bar{x}\} \sim \{\bar{x}\}$. Since $\{\bar{x}, \bar{y}\} \sim \{\bar{y}, \bar{x}\}$, it follows by transitivity that $\{\bar{x}\} \sim \{\bar{y}\}$.

$(\Leftarrow)$ Suppose $\{\bar{x}\} \sim \{\bar{y}\}$. Then the immediately previous argument implies $\{\overline{(\bar{x})}\} \sim \{\overline{(\bar{y})}\}$. Since $\{\overline{(\bar{x})}\} \sim \{x\}$ and $\{\overline{(\bar{y})}\} \sim \{y\}$ by Proposition 3, transitivity of $\sim$ implies $\{x\} \sim \{y\}$. ∎

Definition 3: We say that $x \to y$ if and only if $\{x, \bar{y}\} \notin \wp$.

Proposition 8: $x \to x$.

Proof: Since by Definition 2 we have that $\{x, \bar{x}\} \notin \wp$, by Definition 3 it follows that $x \to x$. ∎

Proposition 9: $\{x\} \sim \{y\}$ if and only if $x \to y$ and $y \to x$.

Proof: $(\Rightarrow)$ Suppose $\{x\} \sim \{y\}$. By the definition of $\sim$ it follows that $\{x, \bar{y}\} \sim \{y, \bar{y}\}$. By Definition 2 we have $\{y, \bar{y}\} \notin \wp$, thus by the definition of $\sim$ it follows that $\{x, \bar{y}\} \notin \wp$ and thus $x \to y$. Since $\{x\} \sim \{y\}$, it follows by the definition of $\sim$ that $\{x, \bar{x}\} \sim \{y, \bar{x}\}$. By Definition 2 we have $\{x, \bar{x}\} \notin \wp$, thus by the definition of $\sim$ it follows that $\{y, \bar{x}\} \notin \wp$ and thus $y \to x$.

$(\Leftarrow)$ Suppose $x \to y$ and $y \to x$. Then by Definition 3, $\{x, \bar{y}\} \notin \wp$ and $\{y, \bar{x}\} \notin \wp$, hence by Proposition 6, we have $\{x\} \sim \{y\}$. ∎

Proposition 10: If $x \to y$ and $y \to z$, then $x \to z$.

Proof: Suppose $x \to y$ and $y \to z$. Since $x \to y$, we have by Definition 3 that $\{x, \bar{y}\} \notin \wp$, i.e. $\{\bar{y}, x\} \notin \wp$, which by Definition 2 implies $\{y, x\} \sim \{x\}$. The definition of $\sim$ implies $\{y, x, \bar{z}\} \sim \{x, \bar{z}\}$. But $y \to z$, hence $\{y, \bar{z}\} \notin \wp$ by Definition 3, and then immediately $\{y,x,\bar{z}\} \notin \wp$ since $\{y,\bar{z}\} \subseteq \{y,x,\bar{z}\}$. Since $\{y,x,\bar{z}\} \sim \{x,\bar{z}\}$, this implies $\{x,\bar{z}\} \notin \wp$. But this means $x \to z$. ∎

<u>Proposition 11</u>: If $t \to x$ and $t \to y$, if $\{x,y\} \sim \{z\}$ then $t \to z$.

<u>Proof</u>: Suppose $t \to x$ and $t \to y$, and $\{x,y\} \sim \{z\}$. Then $\{t,\bar{x}\} \notin \wp$ and $\{t,\bar{y}\} \notin \wp$ by Definition 3. By Definition 2, $\{t,\bar{x}\} \notin \wp$ implies $\{t,x\} \sim \{t\}$, and $\{t,\bar{y}\} \notin \wp$ implies $\{t,y\} \sim \{t\}$, so $\{x,y,t\} \sim \{x,t\} \sim \{t\}$, and thus by transitivity we have $\{x,y,t\} \sim \{t\}$. By definition of equivalence, this means $\{t,\bar{z}\} \sim \{x,y,t,\bar{z}\}$. Since $\{x,y\} \sim \{z\}$, we have $\{x,y,t,\bar{z}\} \sim \{z,t,\bar{z}\} \supseteq \{\bar{z},z\} \notin \wp$, hence $\{x,y,t,\bar{z}\} \notin \wp$, thus since $\{t,\bar{z}\} \sim \{x,y,t,\bar{z}\}$, we have, again by definition of equivalence, that $\{t,\bar{z}\} \notin \wp$ and thus $t \to z$, as desired. ∎

<u>Proposition 12</u>: If $x \to t$ and $y \to t$, and $\{\bar{x},\bar{y}\} \sim \{\bar{z}\}$ then $z \to t$.

<u>Proof</u>: Suppose $x \to t$ and $y \to t$, and $\{\bar{x},\bar{y}\} \sim \{\bar{z}\}$. Then by Definition 3 we have $\{x,\bar{t}\} \notin \wp$ and $\{y,\bar{t}\} \notin \wp$. By Definition 2, $\{x,\bar{t}\} \notin \wp$ implies $\{\bar{x},\bar{t}\} \sim \{\bar{t}\}$, and $\{y,\bar{t}\} \notin \wp$ implies $\{\bar{y},\bar{t}\} \sim \{\bar{t}\}$. Since $\{\bar{x},\bar{t}\} \sim \{\bar{t}\}$, by definition of equivalence we have $\{\bar{x},\bar{y},\bar{t}\} \sim \{\bar{t},\bar{y}\}$. Since $\{\bar{t},\bar{y}\} \sim \{\bar{t}\}$, this implies $\{\bar{x},\bar{y},\bar{t}\} \sim \{\bar{t}\}$. By definition of equivalence, this gives $\{z,\bar{t}\} \sim \{z,\bar{x},\bar{y},\bar{t}\}$. Since $\{\bar{x},\bar{y}\} \sim \{\bar{z}\}$, again by definition of equivalence we get $\{z,\bar{x},\bar{y},\bar{t}\} \sim \{z,\bar{z},\bar{t}\}$. Combining these results yields $\{z,\bar{t}\} \sim \{z,\bar{z},\bar{t}\}$. Since $\{z,\bar{z},\bar{t}\} \supseteq \{z,\bar{z}\} \notin \wp$, it follows that $\{z,\bar{z},\bar{t}\} \notin \wp$. Since $\{z,\bar{t}\} \sim \{z,\bar{z},\bar{t}\}$, we get that $\{z,\bar{t}\} \notin \wp$ and thus $z \to t$. Thus in a sense $z$ is a least upper bound of $x$ and $y$. ∎

<u>Definition 4:</u> If there exists some $z$ such that $\{\bar{x},\bar{y}\} \sim \{z\}$, we define the <u>union</u> $x \vee y = \bar{z}$.

<u>Proposition 13</u>: $x \to y$ if and only if $\bar{y} \to \bar{x}$.

Proof: ($\Rightarrow$) Suppose $x \to y$. Then by Definition 3, $\{x, \overline{y}\} \notin \wp$, hence $\{\overline{y}, x\} \notin \wp$, and $\{\overline{y}, \overline{(\overline{x})}\} \notin \wp$ by Proposition 7 and the definition of equivalence. But this means $\overline{y} \to \overline{x}$.

($\Leftarrow$) Suppose $\overline{y} \to \overline{x}$. Then by Definition 3, $\{\overline{y}, \overline{(\overline{x})}\} \notin \wp$, hence $\{\overline{(\overline{x})}, \overline{y}\} \notin \wp$ and $\{x, \overline{y}\} \notin \wp$ by Proposition 7 and the definition of equivalence. But this means $x \to y$. ∎

<u>Proposition 14</u>: If $x \to y$, and there exist $z, w$ such that $\{t, x\} \sim \{z\}$ and $\{t, y\} \sim \{w\}$, then $z \to w$.

<u>Proof</u>: Since $x \to y$, by Definition 3 we have $\{x, \overline{y}\} \notin \wp$, and thus $\{y, x\} \sim x$ by Definition 2. The definition of equivalence implies $\{t, y, x\} \sim \{t, x\}$. Since $\{y, x\} \sim x$, equivalence implies $\{t, x, \overline{w}\} \sim \{t, y, x, \overline{w}\}$. Since $\{t, y\} \sim \{w\}$ and thus equivalence implies $\{t, y, \overline{w}\} \sim \{w, \overline{w}\}$; by Definition 2, $\{w, \overline{w}\} \notin \wp$ so we get $\{t, y, \overline{w}\} \notin \wp$. Since $\{t, y, x, \overline{w}\} \supseteq \{t, y, \overline{w}\}$, we get $\{t, y, x, \overline{w}\} \notin \wp$. Since $\{t, x, \overline{w}\} \sim \{t, y, x, \overline{w}\}$, this implies $\{t, x, \overline{w}\} \notin \wp$. Since $\{t, x\} \sim \{z\}$, and thus by equivalence $\{t, x, \overline{w}\} \sim \{z, \overline{w}\}$, hence $\{t, x, \overline{w}\} \notin \wp$ implies $\{z, \overline{w}\} \notin \wp$. But this just means that $z \to w$, as desired. ∎

## 2. REFERENCES

[1] Boole, George (1848), The Calculus of Logic, *The Cambridge and Dublin Mathematical Journal*, vol. 3.

E-mail address: ksoileau@yahoo.com